\documentclass[preprint, 12pt]{elsarticle}

\usepackage[utf8]{inputenc}
\usepackage{amsmath}
\usepackage{amssymb}
\usepackage{amsthm}
\usepackage{csquotes}
\usepackage{color}
\usepackage[shortlabels]{enumitem}

\theoremstyle{definition}
\newtheorem{theorem}{Theorem}
\newtheorem{lemma}[theorem]{Lemma}
\newtheorem{prop}{Proposition}
\newtheorem{corollary}{Corollary}

\numberwithin{equation}{section}
\usepackage{hyperref}

\journal{Finite Fields and Their Applications}

\begin{document}

\begin{frontmatter}



\title{Prime Power Residue and Linear Coverings of Vector Space over $\mathbb{F}_{q}$}


\author[1]{Bhawesh Mishra}

\affiliation[1]{organization={Department of Mathematics, The Ohio State University},
            addressline={$231$ W $18^{th}$ Ave}, 
            city={Columbus},
            postcode={43210}, 
            state={OH},
            country={USA}}

\begin{abstract}
Let $q$ be an odd prime and $B = \{b_{j}\}_{j=1}^{l}$ be a finite set of nonzero integers that does not contain a perfect $q^{th}$ power. We show that $B$ has a $q^{th}$ power modulo every prime $p \neq q$ and not dividing $\prod_{b\in B} b$ if and only if $B$ corrresponds to a linear hyperplane covering of $\mathbb{F}_{q}^{k}$. Here, $k$ is the number of distinct prime factors of the $q$-free part of elements of $B$. Consequently: $(i)$ a set $B \subset\mathbb{Z}\setminus\{0\}$ with cardinality less than $q+1$ cannot have a $q^{th}$ power modulo almost every prime unless it contains a perfect $q^{th}$ power and $(ii)$ For every set $B = \{b_{j}\}_{j=1}^{l} \subset\mathbb{Z}\setminus\{0\}$ and for every $\big(c_{j}\big)_{j=1}^{l} \in\Big(\mathbb{F}_{q}\setminus\{0\}\Big)^{l}$ the set $B$ contains a $q^{th}$ power modulo every prime $p \neq q$ and not dividing $\prod_{j=1}^{l}$ if and only if the set $\{b_{j}^{c_{j}}\}_{j=1}^{l}$ does so.
\end{abstract}

\end{frontmatter}



\section{Introduction}

\subsection{History and Motivation.}
Our starting point is the following classical result on sets containing a quadratic residue modulo almost every prime, which was first obtained by Fried in \cite{Fr}, and later rediscovered by Filaseta and Richman in \cite{FR}. 
\begin{prop}
Let $a_{1}, a_{2}, \ldots , a_{l}$ be finitely many distinct nonzero integers. Then the following three conditions are equivalent:
\begin{enumerate}
    
    \item The set $\{a_{1}, a_{2}, \ldots, a_{l}\}$ contains a square modulo almost every prime.\vspace{1mm} 
    
    \item For every odd prime $p \nmid \prod_{j=1}^{l} a_{j}$, the set $\{a_{1}, a_{2}, \ldots, a_{l}\}$ contains a square modulo $p$.\vspace{1mm}
    
    \item There exists $T \subseteq \{1, 2, \ldots , l\}$ of odd cardinality such that $\prod_{j\in T} a_{j} $ is a perfect square.
    
\end{enumerate}
\label{prop1}
\end{prop}
A. Schinzel and M. Skalba substantially generalized Proposition \ref{prop1} by obtaining the necessary and sufficient conditions for finite subsets of a number field $K$ to contain a $n^{th}$ power ($n \geq 2$) modulo almost every prime \cite[Theorems 1, 2]{SS}. However, these conditions were quite complex for a general $n^{th}$ power. In a followup paper, M. Skalba further simplified these conditions in the case $n = q$ is a prime \cite{Sk}. When $K = \mathbb{Q}$, Skalba's result is as follows. 

\begin{prop}
Let $q$ be a prime and $B = \{ b_{1}, b_{2}, \ldots, b_{l}\}$ be a set of finitely many distinct non-zero integers. Then the following conditions are equivalent:
\begin{enumerate}

    \item The set $B$ contains a $q^{th}$ power modulo $p$ for almost every prime $p$. \vspace{1mm}
    
    \item For every prime $p \neq q$ and $p \nmid \prod_{j=1}^{l} b_{j}$, the set $B$ contains a $q^{th}$ power modulo almost every prime. \vspace{1mm}

    \item For each sequence of integers $\{c_{j}\}_{j=1}^{l}$, there exists a sequence of integers $\{f_{j}\}_{j=1}^{l}$ such that \begin{equation}
\sum_{j=1}^{l} f_{j} \not\equiv 0 (\text{ mod } q) \text{ and } \prod_{j=1}^{l} b_{j}^{c_{j}f_{j}} = d^{q} \text{ for some integer } d. \label{skalbalemmaequation}
\end{equation}
\end{enumerate}
\label{Skalbalemma}
\end{prop}
The main result in \cite{SS} deals with a more general problem than sets containing certain power modulo almost every prime, and hence only states the equivalence between condition $(1)$ and $(3)$. However, it follows from results in \cite{SS} that both conditions $(1)$ and $(3)$ are in fact equivalent to the condition $(2)$ in Proposition \ref{Skalbalemma}. Recall that if a number field $K$ contains the complex $q^{th}$ root of unity $\zeta_{q}$, then for every prime ideal $\mathfrak{p}$ of $K$ coprime to $q$ and every $\mathfrak{p}$-adic unit $\alpha \in K$, we define the $q^{th}$ power residue symbol $\big(\alpha|\mathfrak{p}\big)_{q}$ to be the unique $q^{th}$ root of unity $\zeta_{q}^{j}$ such that 
\[\alpha^{\frac{\text{Norm}(\mathfrak{p})-1}{q}} \equiv \zeta_{q}^{j} \hspace{1mm} (\text{mod } \mathfrak{p}).\] These conditions derived in \cite{SS} are obtained by working with the power residue symbol. Therefore, we require the exception $p \neq q$ and $p\nmid \prod_{j=1}^{l} b_{j}$ in the condition $2$ of Proposition \ref{Skalbalemma} above. This is to ensure that $p$ is coprime to $q$ and the integers $b_{1}, b_{2}, \ldots, b_{l}$ are $p$-adic units modulo $p$, so that all the power residue symbols $\big(b_{j}|p\big)_{q}$, as $1 \leq j \leq l$, are well-defined.

One can easily see that the condition $(3)$ in \ref{prop1} is in fact equivalent to the condition $(3)$ in \ref{Skalbalemma} when $q=2$. This is because for a sequence $\{f_{j}\}_{j=1}^{l}$ of integers satisfying $\sum_{j=1}^{l} f_{j} \not\equiv 0 \hspace{1mm}(\text{mod } 2)$ is equivalent to the fact that $f_{j} \equiv 1 \hspace{1mm} (\text{mod } 2)$ for an odd number of entries $f_{j}$. 

We start the discussion of our main result through an example. Consider the set \begin{equation}
   B = \{p_{1}, p_{2}, p_{1}p_{2}, p_{1}p_{2}^{2}\} \label{introexample}
\end{equation} 
for two distinct primes $p_{1}$ and $p_{2}$ different from $3$. One can easily see that for any prime $p \notin\{3, p_{1}, p_{2}\}$ the set $B$ contains a cubic residue modulo $p$. For instance, if neither $p_{1}$ nor $p_{2}$ is a cube modulo $p$, then two cases arise. The first case is when $p_{1}$ and $p_{2}$ lie in the different coset of $\big(\mathbb{F}_{p}^{\times}\big)$ modulo $\big(\mathbb{F}_{p}^{\times}\big)^{3}$. In this case, $p_{1}p_{2}$ is a cube modulo $p$. In the remaining case, $p_{1}p_{2}^{2}$ is a cube modulo $p$. 

We will show that set $B$ in the example above corresponds to a linear covering (defined in Subsection \ref{covering}) of $\mathbb{F}_{3}^{2}$. More generally, the goal of this article to show that for an odd prime $q$, the property of a finite set $B \subset\mathbb{Z}\setminus\{0\}$ to contain a $q^{th}$ power modulo almost every prime is equivalent to $B$ corresponding to a linear covering of $\mathbb{F}_{q}^{k}$. Here $k$ is the number of distinct primes dividing $q$-free parts of elements of $B$. We also demonstrate in Section \ref{CorEx} that this approach yields interesting new insights into the structure of subsets that contain $q^{th}$ power residue modulo cofinitely many primes. Before we state and prove our main result, we will introduce some notation and elementary results needed to do so. 

\subsection{Adjustment within positive $q^{th}$ power equivalence class.} \label{adjustment}
Unless otherwise stated, $q$ will denote an odd prime in this article from now onwards. We assume that the set $B = \{b_{j}\}_{j=1}^{l} \subset\big(\mathbb{Z}\setminus\{0\}\big)$ does not contain a perfect $q^{th}$ power. Otherwise the set $B$ trivially contains a $q^{th}$ power modulo almost every prime. Since $-1$ is a $q^{th}$ power when $q$ is an odd prime, a negative integer $b$ is a $q^{th}$ power modulo $p$ if and only if $-b$ is a $q^{th}$ power modulo $p$. Therefore, without loss of generality, we can assume that the set $B = \{b_{j}\}_{j=1}^{l}$ only contains positive integers.

Given a positive integer $b$ that is not a perfect $q^{th}$ power and has the unique prime factorization $\prod_{i=1}^{m} p_{i}^{a_{i}}$, the $q$-free part of $b$ can be defined as \[ \text{rad}_{q}(b) := \prod_{i=1}^{m} p_{i}^{a_{i} (\text{mod }q)}.\] Hence, as long as we are concerned with power residues modulo primes $p$, with $p \neq q$ and $p \nmid \prod_{j=1}^{l} b_{j}$, in $B = \{b_{j}\}_{j=1}^{l} \subset\mathbb{Z}\setminus\{0\}$ we can replace the set $B$ by the set \[\Big\{\text{rad}_{q}\big(|b_{j}|\big)\Big\}_{j=1}^{l}.\] For instance, when $q = 3$ one can replace the set $\{24, -104, 54\}$ by $\{3, 13, 2\}$ because $24 = 2^{3} \cdot 3$, $104 = 2^{3} \cdot 13$, and $54 = 3^{3} \cdot 2$. 
    
\subsection{Linear Coverings of a Vector Space.} \label{covering} Let $V$ be a vector space of dimension $k \geq 2$ over a field $K$. A linear covering of $V$ is a collection of proper subspaces $\{W_{i}\}_{i\in I}$ such that $V = \bigcup_{i \in I} W_{i}$. The linear covering number of a vector space $V$, denoted by $\#$ LC$(V)$, is the minimum cardinality of a linear covering of $V$. We will use the following fact about $\#$ LC$(V)$, which is the part of the main result proved in \cite{Clark}.
\begin{prop}
For every $\mathbb{F}_{q}$ vector space $V$ of dimension $\geq 2$, we have that \#LC(V) = q + 1. \label{ILC}
\end{prop}
We will also use a well-known result that characterizes linear dependence of a functional on a vector space in terms of containment of null-spaces. More specifically, we will use the following result as given in \cite[pp. 110, Theorem 20]{HK}.
\begin{prop}
    Let $g, f_{1}, f_{2}, \ldots, f_{r}$ be linear functionals on a vector space $V$ with respective null spaces $N, N_{1}, N_{2}, \ldots, N_{r}$. Then, $g$ is a linear combination of $f_{1}, f_{2}, \ldots, f_{r}$ if and only if $\bigcap_{j=1}^{r} N_{j} \subseteq N$. \label{lindep}
\end{prop}
The final ingredient we need to state our result is identification of a finite subset of integers with a certain set of hyperplanes in a $\mathbb{F}_{q}$ vector space. 

\subsection{Finite subsets of integers and hyperplanes in a vector space.}
Let $p_{1}, p_{2}, \ldots, p_{k}$ be finitely many primes that divide any of the element in the set $\{\text{rad}_{q}(|b_{j}|) : 1 \leq j \leq l\}$. Let $\nu_{ij} \geq 0$ such that $p_{i}^{\nu_{ij}} \mid\mid \text{rad}_{q}(|b_{j}|)$ for every $1 \leq i \leq k$ and every $1 \leq j \leq l$. Consider the equations of hyperplanes in $\mathbb{F}_{q}^{k}$ defined as \[\sum_{i=1}^{k} \nu_{ij} x_{i} = 0 \text{ for } 1 \leq j \leq l,\] in the variables $x_{1}, x_{2}, \ldots, x_{k}$. We call the set \[\Big\{ \sum_{i=1}^{k} \nu_{ij} x_{i} = 0 : 1 \leq j \leq l \Big\},\] the set of hyperplanes in $\mathbb{F}_{q}^{k}$ associated with $B = \{b_{1}, b_{2}, \ldots, b_{l}\}$. 

We make the assumption that $k \geq 2$ in the above theorem since one can easily show that for any prime $p_{1}$, there exist infinitely many primes $p$ such that the set $B = \{p_{1}, p_{1}^{2}, \ldots, p_{1}^{q-1}\}$ does not contain a $q^{th}$ power modulo $p$. For instance, we note that the Galois extension $\mathbb{Q}\big(p_{1}^{1/q}, \zeta_{q}\big)\big/\mathbb{Q}$ is the splitting field of the irreducible polynomial $f(x) = x^{q} - p_{1}$. Here $\zeta_{q}$ is a primitive $q^{th}$ root of unity. The Galois group $G$ of this extension is semi-direct product of $\big(\mathbb{Z}/q\mathbb{Z}\big)$  and $\big(\mathbb{Z}/q\mathbb{Z}\big)^{\times}$. The automorphism group $H$ of the intermediate extension $\mathbb{Q}\big(p_{1}^{1/q}, \zeta_{q}\big)\big/\mathbb{Q}\big(p_{1}^{1/q}\big)$ is isomorphic to $\big(\mathbb{Z}/q\mathbb{Z}\big)^{\times}$ and forms a conjugacy class $\mathcal{C}$ in $G$. Therefore, an application of the Chebotarev's density theorem (for instance, as given in \cite[pp. 58, Theorem 5.6]{FrJar}) gives that the set of primes modulo which the polynomial $f$ is still irreducible has density equal to $\frac{|\mathcal{C}|}{|G|} = \frac{q-1}{q(q-1)} = \frac{1}{q}$. In other words, the set of primes modulo which $p_{1}$ is not a $q^{th}$ power has density $\frac{1}{q}$. Our main result is the following theorem. 
\begin{theorem}
Let $q$ be an odd prime and let $B = \{b_{1}, b_{2}, \ldots, b_{l}\}\subset \mathbb{Z}\setminus\{0\}$ not containing a perfect $q^{th}$ power. Let $p_{1}, p_{2}, \ldots, p_{k}$ $(k \geq 2)$ be all the primes dividing elements of $\{\text{rad}_{q}(|b_{j}|) : 1 \leq j \leq l\}$ such that $p_{i}^{\nu_{ij}} \mid\mid \text{rad}_{q}(|b_{j}|)$ for every $1 \leq i \leq k$ and every $1 \leq j \leq l$. Then, the following three conditions are equivalent:
\begin{enumerate}
    \item The set $B$ contains a $q^{th}$ power modulo almost every prime. \vspace{1mm}
    
    \item For every prime $p\neq q$ and not dividing $\prod_{j=1}^{l} b_{j}$, the set $B$ contains a $q^{th}$ power modulo $p$. 
    
    \item The set of hyperplanes \[\Big\{ \sum_{i=1}^{k} \nu_{ij} x_{i} = 0 : 1 \leq j \leq l \Big\}\] associated with $B$, contains a linear covering of $\mathbb{F}_{q}^{k}$.
\end{enumerate}\label{Mainresult2}
\end{theorem}

\section{Proof of Theorem \ref{Mainresult2}}\label{restatement}
As mentioned in the Subsection \ref{adjustment}, without loss of generality, we can assume that the elements of $B$ are positive and $q$-free. Therefore, we will replace $\text{rad}_{q}\big(|b_{j}|\big)$ by $b_{j}$ in the statement of Theorem \ref{Mainresult2} from now onwards. Since for every $1 \leq i \leq k$ and for every $1 \leq j \leq l$, $p_{i}^{\nu_{ij}} \mid\mid b_{j}$, we have $b_{j} = \prod_{i=1}^{k} p_{i}^{\nu_{ij}}$ for every $1 \leq j \leq l$. Therefore, the condition \[\prod_{j=1}^{l}b_{j}^{c_{j}f_{j}} = d^{q} \] in equation \ref{skalbalemmaequation} can be rewritten as 
\begin{multline}
    \prod_{j=1}^{l} \Big( \prod_{i=1}^{k} p_{i}^{\nu_{ij}} \Big)^{c_{j}f_{j}} = d^{q}  \Leftrightarrow \prod_{j=1}^{l}    \Big( \prod_{i=1}^{k} p_{i}^{\nu_{ij}c_{j}f_{j}} \Big) = d^{q}\\
    \Leftrightarrow \sum_{j=1}^{l} \nu_{ij}c_{j}f_{j} \equiv 0 \hspace{1mm} (\text{mod } q) \text{ for every } 1 \leq i \leq k. \label{restatement1}
\end{multline}
Consequently, it suffices to consider $c_{j}, f_{j} \in \mathbb{F}_{q}$ in the condition $(2)$ of Proposition \ref{Skalbalemma}. We also note that if $c_{j_{0}} = 0$ for some $1 \leq j_{0} \leq l$, then the equation \eqref{restatement1}, and hence condition $(2)$ in Proposition \ref{Skalbalemma}, can be readily solved by taking $f_{j} = \begin{cases} 0 \text{ for } j \neq j_{0} \\ 1 \text{ for } j = j_{0}\end{cases}$. Therefore, we only need to consider $(c_{j})_{j=1}^{l} \in\big(\mathbb{F}_{q}\setminus\{0\}\big)^{l}$. 

The system of congruences \[ \sum_{j=1}^{l} \nu_{ij} c_{j} f_{j} \equiv 0 \hspace{1mm} (\text{mod } q), \hspace{1mm} 1 \leq i \leq k\] can be written as 
\begin{equation}
M\Big((c_{j})_{j=1}^{l} \Big) \times F = 0 \label{MF} 
\end{equation} 
over the field $\mathbb{F}_{q}$, where $M\Big((c_{j})_{j=1}^{l} \Big)$ is $(k \times l)$ matrix whose $ij$-th entry is given by$M_{ij} := \nu_{ij} c_{j}$ and $F = \begin{bmatrix} f_{1} \\ f_{2} \\ \vdots \\ f_{l} \end{bmatrix}$. Now, we prove a lemma that replaces the condition $3$ in Proposition \ref{Skalbalemma} by another equivalent statement. 

\begin{lemma}
Using the notation defined above in Section \ref{restatement}, the following three statements are equivalent.
\begin{enumerate}[A.]

    \item The condition $3$ in Proposition \ref{Skalbalemma}.

    \item For every $(c_{j})_{j=1}^{l} \in\big(\mathbb{F}_{q}\setminus\{0\}\big)^{l}$, there exists $(f_{j})_{j=1}^{l} \in\big(\mathbb{F}_{q}\setminus\{0\}\big)^{l}$ such that $\sum_{j=1}^{l} \nu_{ij}c_{j}f_{j} \equiv 0 \hspace{1mm} (\text{mod } q)$ for every $1 \leq i \leq k$ and $\sum_{j=1}^{l} f_{j} \not\equiv 0 \hspace{1mm} (\text{mod } q)$. 
    
    \item For every $(c_{j})_{j=1}^{l} \in\big(\mathbb{F}_{q}\setminus\{0\}\big)^{l}$, the row space of the matrix $M\Big((c_{j})_{j=1}^{l} \Big)$ does not contain $\begin{bmatrix} 1 & 1 & \ldots & 1 \end{bmatrix}$.
\end{enumerate}
\end{lemma}

\begin{proof}
Let $(c_{j})_{j=1}^{l} \in\big(\mathbb{F}_{q}\setminus\{0\}\big)^{l}$. The equivalence between conditions $A$ and $B$ follows from Section \ref{restatement} because $\sum_{j=1}^{l} \nu_{j}c_{j}f_{j} \equiv 0 \hspace{1mm}(\text{mod } q)$ is equivalent to the fact that $F = \begin{bmatrix} f_{1} \\ f_{2} \\ \vdots \\ f_{l} \end{bmatrix}$ lies in the null space of $M\Big((c_{j})_{j=1}^{l} \Big)$ as remarked in \eqref{MF}. Now we prove the equivalence of conditions $B$ and $C$. 

Since any row $R$ of size $1 \times l$ induces a linear functional $v \longmapsto R \times v$ on $\mathbb{F}_{q}^{l}$, using Proposition \ref{lindep} we have that the condition $C$ is equivalent to the fact that \[\text{Null}\Bigg( M\Big((c_{j})_{j=1}^{l} \Big) \Bigg) = \bigcap_{i=1}^{k} \text{Null}(R_{i}) \not\subseteq \text{Null} (\begin{bmatrix} 1 & 1 & \ldots & 1 \end{bmatrix}).\] 

Here, $R_{i}$ is the $i^{th}$ row of the matrix $ M\Big((c_{j})_{j=1}^{l} \Big)$ for $1 \leq i \leq k$. In other words, the condition $C$ is equivalent to the fact that there exists \[F = \begin{bmatrix} f_{1} \\ f_{2} \\ \vdots \\ f_{l} \end{bmatrix} \in \bigcap_{i=1}^{k} \text{Null}(R_{i}) = \text{Null}\Bigg( M\Big((c_{j})_{j=1}^{l} \Big) \Bigg)\] such that $\sum_{j=1}^{l} f_{j} \not\equiv 0 \hspace{1mm} (\text{mod } q)$.
\end{proof}

Therefore, as a result of Proposition \ref{Skalbalemma}, to prove the Theorem \ref{Mainresult2} we need to prove equivalence of following two statements.
\begin{enumerate}[(I)]
    \item For every $(c_{j})_{j=1}^{l} \in \Big( \big(\mathbb{F}_{q}\big) \setminus\{0\}\Big)^{l}$ the linear span of rows of $M\Big((c_{j})_{j=1}^{l} \Big)$ does not contain the vector $\begin{bmatrix} 1 & 1 & \ldots & 1 \end{bmatrix}$. 
    
    \item The set of hyperplanes \[\Big\{ \sum_{i=1}^{k} \nu_{ij} x_{i} = 0 : 1 \leq j \leq l \Big\}\] associated with $B$, contains a linear covering of $\mathbb{F}_{q}^{k}$.
\end{enumerate}

\subsection{Proof of (II) implies (I)} Assume that the set of hyperplanes \[\Big\{ \sum_{i=1}^{k} \nu_{ij} x_{i} = 0 : 1 \leq j \leq l \Big\}\] associated with $B$, contains a linear covering of $\mathbb{F}_{q}^{k}$. We are left to show that no $\big(\mathbb{F}_{q}\setminus\{0\}\big)$-linear combination of rows of $M\Big((c_{j})_{j=1}^{l}\big)$ is equal to  $\begin{bmatrix} 1 & 1 & \ldots & 1 \end{bmatrix}$.

Let $(c_{j})_{j=1}^{l} \in \big(\mathbb{F}_{q}\setminus\{0\}\big)^{l}$ and let $R_{1}, R_{2}, \ldots, R_{k}$ be the rows of the matrix $M\big((c_{j})_{j=1}^{l}\big)$. In other words,
\[R_{i} := \begin{bmatrix} c_{1}\nu_{i1} & c_{2}\nu_{i2} & \ldots & c_{l}\nu_{il} \end{bmatrix} \text{ for every } 1 \leq i \leq k.\]
Let $d_{1}, d_{2}, \ldots, d_{k} \in\mathbb{F}_{q}$ and consider the linear combination $d_{1}R_{1} + d_{2}R_{2} + \ldots + d_{k} R_{k}$ of rows, which is equal to 
\[\begin{bmatrix}c_{1}\big(\sum_{i=1}^{k} d_{i}\nu_{i1}\big) & c_{2}\big(\sum_{i=1}^{k} d_{i}\nu_{i2}\big) & \ldots & c_{1}\big(\sum_{i=1}^{k} d_{i}\nu_{il}\big)  \end{bmatrix}.\]
Since, the set of hyperplanes $\Big\{ \sum_{i=1}^{k} \nu_{ij} x_{i} = 0 : 1 \leq j \leq l \Big\}$ contains a linear covering of $\mathbb{F}_{q}^{k}$ and $(d_{i})_{i=1}^{k} \in\mathbb{F}_{q}^{k}$, there exists $1 \leq j_{0} \leq l$ such that $\sum_{i=1}^{k} \nu_{ij_{0}} d_{i} = 0$. In other words, one of the entries of $d_{1}R_{1} + d_{2} R_{2} + \ldots + d_{k} R_{k}$ is equal to zero, and hence cannot be equal to the row matrix  $\begin{bmatrix} 1 & 1 & \ldots & 1 \end{bmatrix}$. Since $(c_{j})_{j=1}^{l} \in \big(\mathbb{F}_{q}\setminus\{0\}\big)^{l}$ and $d_{1}, d_{2}, \ldots, d_{k} \in\mathbb{F}_{q}$ were arbitrary, the proof of this direction is complete. $\square$

\subsection{Proof of (I) implies (II)} Let us assume that the set of hyperplanes \[\Big\{ \sum_{i=1}^{k} \nu_{ij} x_{i} = 0 : 1 \leq j \leq l \Big\}\] does not contain a linear covering of $\mathbb{F}_{q}^{k}$. We shall show that for some $(c_{j})_{j=1}^{l} \in\big(\mathbb{F}_{q}\setminus\{0\}\big)^{l}$, some linear combination of the rows of the $(k \times l)$ matrix $M\big((c_{j})_{j=1}^{l}\big)$ is equal to $\begin{bmatrix} 1 &  1 & 1 \ldots & 1\end{bmatrix}$.

Since the set $\Big\{ \sum_{i=1}^{k} \nu_{ij} x_{i} = 0 : 1 \leq j \leq l \Big\}$ does not contain a linear covering of $\mathbb{F}_{q}^{k}$, there exists $(d_{1}, d_{2}, \ldots, d_{k})\in\mathbb{F}_{q}^{k}$ such that \[\sum_{i=1}^{k} \nu_{ij}d_{i} \neq 0 \text{ for any } 1 \leq j \leq l.\] Define \[c_{j} := \big( \sum_{i=1}^{k} \nu_{ij}d_{i} \big)^{-1} \in\mathbb{F}_{q}\setminus\{0\} \text{ for every } 1 \leq j \leq l.\] Consider the matrix $M\big((c_{j})_{j=1}^{l}\big)$ associated with $(c_{j})_{j=1}^{l}$, and let $R_{1}, R_{2}, \ldots, R_{k}$ be its rows, i.e., 
\[R_{i} := \begin{bmatrix} c_{1}\nu_{i1} & c_{2}\nu_{i2} & \ldots & c_{l}\nu_{il} \end{bmatrix} \text{ for every } 1 \leq i \leq k.\]
The linear combination $d_{1}R_{1}+d_{2}R_{2}+\ldots+d_{k}R_{k}$ of rows, as in previous subsection, is equal to 
\[\begin{bmatrix}c_{1}\big(\sum_{i=1}^{k} d_{i}\nu_{i1}\big) & c_{2}\big(\sum_{i=1}^{k} d_{i}\nu_{i2}\big) & \ldots & c_{1}\big(\sum_{i=1}^{k} d_{i}\nu_{il}\big)  \end{bmatrix} = \begin{bmatrix} 1 & 1 & \ldots & 1\end{bmatrix}.\] The last equality follows from definition of $c_{j}$. $\square$ 

\section{Corollaries and Examples.}\label{CorEx}
\subsection{Corollaries.}
Our first corollary puts a lower bound on cardinality of subsets of integers that contain a $q^{th}$ power modulo almost every prime but do not contain a perfect $q^{th}$ power. 
\begin{corollary}
If a set $B \subset\mathbb{Z}$ with $|B| < q+1$ satisfies the conditions in Theorem \ref{Mainresult2}, then $B$ contains a perfect $q^{th}$ power. More specifically, a set of non-zero integers with cardinality less than $q+1$ contains a $q^{th}$ power modulo almost every prime if and only if $B$ contains a perfect $q^{th}$ power. \label{cor1}
\end{corollary}

\begin{proof}
Assuming $|B| = l$ and that $B$ does not contain a perfect $q^{th}$ power, the Theorem \ref{Mainresult2} implies that the set of $l$ hyperplanes associated with $B$ must contain a linear covering of $\mathbb{F}_{q}^{k}$. However, Proposition \ref{ILC} implies that any such linear covering must contain at least $q+1$ elements, implying that $|B| = l \geq q + 1$. 
\end{proof}
Theorem \ref{Mainresult2} also implies that the property whether a given finite set $B$ of integers contains a $q^{th}$ power modulo almost every prime is invariant under exponentiation by elements of $\mathbb{F}_{q}\setminus\{0\}$ in the following sense. 

\begin{corollary}
Let $B = \{b_{1}, b_{2}, \ldots, b_{l}\} \subset \mathbb{Z}$ , $(a_{j})_{j=1}^{l} \in \big( \mathbb{F}_{q}\setminus\{0\} \big)^{l}$ and \[B^{\prime} := \big\{ b_{j}^{a_{j}} : 1 \leq j \leq l \big\}.\] Then, $B$ satisfies the conditions in Theorem \ref{Mainresult2} if and only if $B^{\prime}$ satisfies the conditions in Theorem \ref{Mainresult2}. \label{cor2}
\end{corollary}

\begin{proof}
Without loss of generality, we can assume that $B$ does not contain a perfect $q^{th}$ power. Otherwise, $B^{\prime}$ also contains a perfect $q^{th}$ power and both $B$ and $B^{\prime}$ contain a $q^{th}$ power modulo almost every prime. 

Let $p_{1}, p_{2}, \ldots, p_{k}$ be all the primes that divide an element of $\big\{ \text{rad}_{q}(|b_{j}|)\big\}_{1 \leq j \leq l}$. For every $1 \leq j \leq l$ and every $1 \leq i \leq k$, let $\nu_{ij} \geq 0$ such that $p_{i}^{\nu_{ij}} \mid\mid b_{j}$. Then, the set of hyperplanes in $\mathbb{F}_{q}^{k}$ associated with the set $B$ is equal to \[H:= \big\{ \sum_{i=1}^{k} \nu_{ij} x_{i} = 0 \text{ for } 1 \leq j \leq l. \big\}\] Since, $p_{i}^{\nu_{ij}} \mid\mid b_{j}$, we have that $p_{i}^{a_{j}\nu_{ij}} \mid\mid b_{j}^{a_{j}}$. Therefore, the set of hyperplanes in $\mathbb{F}_{q}^{k}$ associated with the set $B^{\prime}$ is equal to 
\[H^{\prime}:= \big\{ \sum_{i=1}^{k} a_{j}\nu_{ij} x_{i} = 0 \text{ for } 1 \leq j \leq l \big\} = \big\{ \sum_{i=1}^{k} \nu_{ij} x_{i} = 0 \text{ for } 1 \leq j \leq l \big\} = H.\]
Since the set of hyperplanes in $\mathbb{F}_{q}^{k}$ associated with $B$ is equal to the set of hyperplanes in $\mathbb{F}_{q}^{k}$ associated with $B^{\prime}$, the Theorem \ref{Mainresult2} implies the corollary.
\end{proof}

\subsection{Examples.}As an application of Theorem \ref{Mainresult2}, we now present some illustrative examples of sets that contain a $q^{th}$ power modulo almost every prime for $q = 3$ and $q = 5$.
\begin{enumerate}
    \item Now we will use Theorem $1$ to justify that the set \[B = \{p_{1}, p_{2}, p_{1}p_{2}, p_{1}p_{2}^{2}\}\] as in Example \eqref{introexample} contains a cube modulo every prime $p$ such that $p \nmid 3p_{1}p_{2}$. Here $p_{1}, p_{2}$ be two distinct primes different from $3$.
    
    Note that the set of hyperplanes \[H := \big\{ x_{1} = 0, x_{2} = 0, x_{1} + x_{2} = 0, x_{1} + 2x_{2} = 0 \big\}\] forms a linear covering of $\mathbb{F}_{3}^{2}$, because every $(c_{1}, c_{2}) \in\mathbb{F}_{3}^{2}$ satisfies at least one of the equations in $H$. Therefore, Theorem \ref{Mainresult2} implies that the set \[B = \{p_{1}, p_{2}, p_{1}p_{2}, p_{1}p_{2}^{2}\}\] contains a cube modulo almost every prime. \vspace{1mm}
    
    \item Using $(1)$ and Corollary \ref{cor2}, one could obtain $|\big(\mathbb{F}_{3}\setminus\{0\}\big)^{4}|-1 = 15$ more examples of sets \[\Big\{ \{p_{1}^{a_{1}}, p_{2}^{a_{2}}, p_{1}^{a_{3}}p_{2}^{a_{3}}, p_{1}^{a_{4}}p_{2}^{2a_{4}}\} : (a_{j})_{j=1}^{4} \in\big(\mathbb{F}_{3}\setminus\{0\}\big)^{4} \Big\}\] that satisfy the following:
    
    \begin{itemize}
        \item Each of these sets contain a cube modulo almost every prime. 
        
        \item Only prime divisors of elements of these sets are $p_{1}$ and $p_{2}$. 
        
        \item Each of these sets correspond to the same linear hyperplane covering of $\mathbb{F}_{3}^{2}$ as in $1$.
    \end{itemize}
    
    \item Let $p_{1}, p_{2}, p_{3}$ be three distinct primes different from $5$ and let \[B = \{p_{1}, p_{2}p_{3}, p_{1}p_{2}p_{3}, p_{1}^{2}p_{2}p_{3}, p_{1}^{3}p_{2}p_{3}, p_{1}^{4}p_{2}p_{3} \}.\] Note that
    \begin{multline*}
       \\ \nu_{11} = 1, \nu_{12} = 0, \nu_{13} = 1, \nu_{14} = 2, \nu_{15} = 3, \nu_{16} = 4 \\ \nu_{21} = 0, \nu_{22} = 1, \nu_{23} = 1, \nu_{24} = 1, \nu_{25} = 1, \nu_{26} = 1 \\ \nu_{31} = 0, \nu_{32} = 1, \nu_{33} = 1, \nu_{34} = 1, \nu_{35} = 1, \nu_{36} = 1, \\
    \end{multline*}
    and hence, the set of hyperplanes associated with $B$ is 
    \begin{multline*}
        \Big\{ x_{1} = 0, x_{2}+x_{3} = 0, x_{1}+x_{2}+x_{3} = 0,  2x_{1}+x_{2}+x_{3} = 0, \\ 3x_{1}+x_{2}+x_{3} = 0, 4x_{1}+x_{2}+x_{3} = 0 \Big\},
    \end{multline*}
    which is easily seen to form a linear covering of $\mathbb{F}_{5}^{3}$. This is because every $(x_{j})_{j=1}^{3} \in\mathbb{F}_{5}^{3}$ satisfies one of the equations of hyperplanes above. Therefore, the set $B$ satisfies Theorem \ref{Mainresult2} and contains a fifth power modulo almost every prime.
    
    \item Similar to $(2)$, using $(3)$ and Corollary \ref{cor2}, one could obtain \[|\big(\mathbb{F}_{5}\setminus\{0\}\big)^{6}|-1 = 4^{6} - 1 = 4095 \] more examples of sets that satisfy the following:
    \begin{itemize}
        \item Only prime divisors of elements of these sets are $p_{1}, p_{2}$ and $p_{3}$. 
        
        \item Each of these sets contain a fifth power modulo almost every prime. 
        
        \item Each of these sets correspond to the same linear hyperplane covering of $\mathbb{F}_{5}^{3}$ as in $(3)$.
    \end{itemize}
\end{enumerate}

\section*{Acknowledgement} The present version of this article owes a great deal to both the referees for their invaluable comments and suggestions, which greatly enhanced the exposition in this article.

\end{document}